\documentclass[a4paper,11pt,english]{amsart}

\usepackage{babel,varioref}
\usepackage[ansinew]{inputenc}
\usepackage{amsmath}
\usepackage{amsfonts}
\usepackage{amscd}
\usepackage{enumerate}
\usepackage{graphics}
\usepackage{babel,varioref}
\usepackage{amscd}
\usepackage[pdftex]{graphicx}
\DeclareGraphicsExtensions{.pdf,.png,.jpg}

\input amssym.def
\input amssym.tex

\def\reels{\mathbb{R}}
\def\complexes{\mathbb{C}}

\def\g {\mathfrak}

\newtheorem{theo}{Theorem}
\newtheorem{prop}[theo]{Proposition}

\title[Holonomy direct sum of two faithful representations]{Holonomy representations which are a diagonal direct sum of two faithful representations}
\author[T. Krantz]{Tom Krantz}
\thanks{e-mail: Tom.Krantz@uni.lu\\
Address: Institut Élie Cartan Nancy. Université Henri Poincaré Nancy
1.\\ B.P. 239, F-54506 Vandoeuvre-lès-Nancy Cedex, France.\\
Université du Luxembourg. 162a, avenue de la Faïencerie. L-1511
Luxembourg. G.-D. of Luxembourg.\\ The author is supported by BFR
03/095 of the ministry of research of the G.-D. of Luxembourg.}

\makeindex 

\begin{document}
\begin{abstract}
We study holonomy representations admitting a pair of supplementary
faithful sub-representations. In particular the cases where the
sub-representations are isomorphic respectively dual to each other
are treated. In each case we have a closer look at the
classification in small dimension.
\end{abstract}
\maketitle
\section{Introduction}
The notion of torsion-free connection on the tangent bundle $TM$ of
a smooth manifold $M$ gives rise to the notion of parallel transport
along paths contained in the manifold, and it is well known that if
one considers all closed contractible paths based in the a point $o$
of the manifold $M$, the set of the corresponding parallel
transports forms a Lie subgroup of $Gl(T_o M)$ called the restricted
holonomy group of the connection. To this linear Lie group
corresponds its linear Lie algebra which will be noted $\g g$
throughout this text. Seen as a representation we will refer to it
as the holonomy representation. For a general study of connections
and holonomy and precise definitions we refer for example to
\cite{Kobayashi-Nomizu}.

Holonomy representations $\g g \subset \g{gl}(T_o M)$ have been
completely studied in the Riemannian case through the classical
result of the De Rham-theorem reducing the classification problem to
the one of the classification irreducible metric holonomy
representations. In the pseudo-Riemannian case Wu's generalization
of the de Rham-theorem allows to reduce the classification problem
to the one of weakly-irreducible representations. All irreducible
holonomy representations are by now known, but in the weakly
irreducible case mainly signature $(1,n)$ (see~\cite{LBBAI},
\cite{Ikemakhen2}, \cite{Boubel},
 \cite{leistner1}, \cite{leistner2}, \cite{leistner3},
\cite{leistner4}, \cite{Galaev1}), $(2,n)$ (see~\cite{LBBAI},
~\cite{Ikemakhen3}, \cite{Galaev2}) are explored. For signature
$(n,n)$ see (\cite{LBBAI}). Lionel Bérard Bergery has classified
indecomposable semi-simple non simple pseudo-Riemannian holonomy
representations. These are examples of $V\oplus V^*$ representations
which we examine in this paper generally.

If the connection is not supposed to preserve a metric very little
is known, we have no longer a general result like the De Rham-Wu
theorem. In my thesis(\cite{krantz})(under the direction of Lionel
Bérard Bergery), I studied representations admitting two pairs of
supplementary invariant spaces(see also \cite{LBBTK}). In the
pseudo-Riemannian or symplectic case appear naturally factors of
type $V\oplus V^*$. When the representation is weakly irreducible
and admits a decomposition into two supplementary non trivial
subspaces, it is of this type as shown in the paper. In the general
setting appear factors of type $V\otimes \reels^2$ (with $\reels^2$
being the trivial representation). Here we examine more closely both
types of holonomy representations.

The paper is structured as follows: We start by formulating
Berger-criteria for representations which are a diagonal direct sum
of two faithful representations. We explore then the local geometric
structure of torsion-free connections with this holonomy and
calculate curvature in an adapted basis. Finally we examine closely
holonomy representations of type $V\oplus V^*$ and $V\otimes
\reels^2$ and give classifications in dimension $2\times 2$ (which
was known for the case $V\oplus V^*$ and is new for the case
$V\otimes \reels^2$).

\subsection{Notation}
We will say a representation $\g g \subset \g g \g l(V)$ is {\em
decomposable along the direct sum $V=V_1 \oplus V_2$} if $\g g=\g
g\cap V_1^* \otimes V_1 \oplus \g g\cap V_2^* \otimes V_2$.

\section{Berger-type criteria}
\subsection{The first Berger criterion}

Let $\g g \subset \g g \g l(V)$ be a finite-dimensional
representation.

Note $${\mathcal K}(\g g):= \left\{ R \in (V^* \wedge V^*) \otimes
\g g \; \vrule \; R(x,y)z+R(y,z)x+R(z,x)y=0\mbox{ for }x,y,z\in V
\right\},$$
$$\underline{\g g}:=<R(x,y) \; \vrule \; x,y\in V, R\in
{\mathcal K}(\g g) >.$$

Since the work of Marcel Berger(\cite{berger1}) it is known that:
\begin{prop}
If $\g g \subset \g g \g l(V)$ is a holonomy algebra of a
torsion-free connection then $\underline{\g g}=\g g$.
\end{prop}

Recall also the second criterion of Berger:

Note $${\mathcal K}^1(\g g):=\left\{ D \in V^* \otimes {\mathcal
K}(\g g) \; \vrule \; Dx(y,z)+Dy(z,x)+Dz(x,y)=0\mbox{ for }x,y,z\in
V \right\}.$$

\begin{prop}
If $\g g \subset \g g \g l(V)$ is a holonomy algebra of a
torsion-free connection which is non locally symmetric then
${\mathcal K}^1(\g g)\neq \{0\}$.
\end{prop}

In this text we are interested representations $V$ of a Lie algebra
$\g g$ such that $V=V_1\oplus V_2$ where $V_1$ and $V_2$ are two
faithful representations of $\g g$.

Let's write ${\g g}_i \subset \g g \g l(V_i)$ for the restriction of
the action of $\g g$ to $V_i$.

In the context of these particular representations the first Berger
criterion can be specialized to the following.

Note $$k(\g g):= \left\{ R \in V_1^* \otimes V_2^* \otimes \g g \;
\vrule \; \begin{array}{l}R(x,y')z'=R(x,z')y'\\
R(x,y')t=R(t,y')x\end{array}\mbox{ for }x,t\in V_1, y',z'\in V_2
\right\},$$

We have then:

\begin{prop}
If $\g g \subset \g g \g l(V)$ is a finite dimensional
representation and $V_1$ and $V_2$ are two faithful representations
of $\g g$ such that $V=V_1 \oplus V_2$ then:

$$\underline{\g g}=<R(x,y') \; \vrule \; x\in V_1, y'\in V_2, R\in
k(\g g)>.$$
\end{prop}
\begin{proof}
This follows from the fact that if $R$ is in $K(\g g)$ then
$R(x,y)\in \g g \g l(V_i)$ for $x,y\in V_i$, from which follows
$R(x,y)=0$ by the condition that $V_1$ and $V_2$ are two faithful
representations of $\g g$.
\end{proof}

Note $$k^1(\g g):= \left\{ D \in V^* \otimes k(\g g) \;
\vrule \; \begin{array}{l}Dx(y,z')=Dy(x,z')\\
Dt'(y,z')=Dz'(x,t')\end{array}\mbox{ for }x,y\in V_1, z',t'\in V_2
\right\},$$

Berger's second criterion reads then:

\begin{prop}
Let $\g g \subset \g g \g l(V)$ be a representation and $V_1$ and
$V_2$ two faithful representations of $\g g$ such that $V=V_1 \oplus
V_2$.

If $\g g$ is a holonomy algebra of a torsion-free connection which
is non locally symmetric then $k^1(\g g)\neq \{0\}$.
\end{prop}

\subsection{A general weak criterion}

We write $${\g g}_i^{(1)}:= V_i^* \otimes {\g g}_i \cap
S^2(V_i^*)\otimes V_i,$$

the first prolongation of ${\g g}_i$.

Let $$\underline{\underline{{\g g}_i}}:=<r(x) \; \vrule \; x\in V_i,
\; r\in {\g g}_i^{(1)}>.$$

\begin{prop}
$\g g$ verifies the first Berger criterion implies
$\underline{\underline{{\g g}_i}}={\g g}_i$ for $i=1,2$.
\end{prop}

Weaker than the first Berger criterion this criterion considers only
the ${\g g}_i$ and ignores their coupling.

If $W_i$ is a sub-representation of $V_i$, note ${\g g_i}_{/W_i}
\subset \g g \g l(V_i / W_i)$ the representation obtained from ${\g
g}_i$ by passing to the quotient.

The following result shows some properties of holonomy
representations we are interested in this article in.

\begin{prop}
For $i=1,2$, $\underline{\underline{{\g g}_i}}={\g g}_i$ implies
that for a sub-representation $W_i$ of $V_i$,
$\underline{\underline{{\g g_i}_{/W_i}}}={\g g_i}_{/W_i}$
\end{prop}

\section{Local geometric structure}
\subsection{Basic results}
Assume we have a connection $\nabla$ on a manifold $\mathcal M$
admitting a holonomy representation $\g g \subset \g g \g l(T_o
{\mathcal M})$ (at the point $o\in \mathcal M$) for which $T_o M=V_1
\oplus V_2$, where $V_1$ and $V_2$ are two faithful representations
of $\g g$.

One can transport parallely the direct sum $V_1 \oplus V_2$ to any
point of the manifold. Each $V_i$ gives a distribution on $\mathcal
M$ which is integrable and gives rise to a foliation ${\mathcal
F}_i$ of the manifold with flat leaves. On the leaves ${\mathcal
F}_1(o)$ (resp. ${\mathcal F}_2(o)$) we can choose coordinates
$x_{o}^i$ ($i=1\ldots n_1$) (resp. $y_{o}^j$ ($j=1\ldots n_2$)) for
which $\nabla_{\partial_{x_o^i}}\partial_{x_o^j}=0$, $\forall i,j$,
and $\nabla_{\partial_{y_{o}^i}}\partial_{y_{o}^j}=0$, $\forall
i,j$. The coordinates are in bijection with the basis
$(\partial_{x_o^1},\ldots,\partial_{x_o^{n_1}},\partial_{y_o^1},\ldots,\partial_{y_o^{n_2}})$
in $o$.

The manifold $\mathcal M$ - as a differential manifold - locally a
product and to a point $p$ of the manifold one can (locally)
associate coordinates $(x_o^i(p), i=1\ldots n_1, y_{o}^j(p),
j=1\ldots n_2)$ by considering the coordinates of the intersection
point of ${\mathcal F}_2(p)$ with ${\mathcal F}_1(o)$ (respectively
the ones of the intersection of ${\mathcal F}_1(p)$ with ${\mathcal
F}_2(o)$).

On the other hand one can equip each leaf ${\mathcal F}_1(r)$
containing the point $r\in {\mathcal F}_2(o)$, because it is flat,
with coordinates $x_r^i,i=1\ldots n_1$ such that
$\nabla_{\partial_{x_r^i}}\partial_{x_r^j}=0$, $\forall i,j$ on the
leaf ${\mathcal F}_1(r)$ and such that $\partial_{x_r^i}(r)$ is
obtained by parallely transporting $\partial_{x_o^i}(o)$ along the
leaf ${\mathcal F}_2(o)$. Analogously one defines $y_q^i, i=1\ldots
n_2$ on the leaf ${\mathcal F}_2(q)$.

One has finally as well the following coordinates: $(x^i(p),
i=1\ldots n_1, y^j(p), j=1\ldots n_2)$, with $x^i:=x_o^i$ and
$y^j(p):=y_q^j(p)$.

\begin{prop}\label{p1}
Let $\mathcal M$ be a manifold equipped with a torsion-free
connection$\nabla$ such that the holonomy algebra (in $o$) $\g g$
verifies $T_o M = V_1 \oplus V_2$ with $V_1$ and $V_2$ faithful
representations of $\g g$. For $o$ in $\mathcal M$ let $x_o^i$
(respectively $y_o^j$) designate local flat coordinates on
${\mathcal F}_1(o)$ (resp. ${\mathcal F}_2(o)$) corresponding to a
choice of basis $(b_1)$ of $V_1$ and $(b_2)$ of $V_2$. To $(b_1)$
(resp. $(b_2)$) is associated for any point $r$ of ${\mathcal
F}_2(o)$ (resp. $q$ of ${\mathcal F}_1(o)$) flat coordinates $x_r^i$
(respectively $y_q^j$) of ${\mathcal F}_1(r)$ (resp. ${\mathcal
F}_2(q)$).

Let $p$ be a point of $\mathcal M$ close to $o$ and $q$ (resp. $r$)
the unique intersection point of ${\mathcal F}_1(o)$ with ${\mathcal
F}_2(p)$ respectively ${\mathcal F}_2(o)$ with ${\mathcal F}_1(p)$,
$\gamma_1$ a path contained in ${\mathcal F}_1(o)$ from $o$ to $q$,
$\gamma_2$ a path contained in ${\mathcal F}_2(p)$ from $q$ to $p$,
$\gamma_3$ a path contained in ${\mathcal F}_1(p)$ from $p$ to $r$,
$\gamma_4$ a path contained in ${\mathcal F}_2(o)$ from $r$ to $o$,
$\gamma$ the composition $\gamma_1 \gamma_2 \gamma_3 \gamma_4$. Then
$\tau_\gamma$ depends only on $p$, and one has: $\tau_{\gamma_1
\gamma_2}(\partial_{x_o^i,o})=\partial_{x_o^i,p}$ and
$\tau^{-1}_{\gamma_3
\gamma_4}(\partial_{x_o^i,o})=\partial_{x_r^i,p}$. Similarly
$\tau_{\gamma_1 \gamma_2}(\partial_{y_o^i,o})=\partial_{y_q^i,p}$ et
$\tau^{-1}_{\gamma_3
\gamma_4}(\partial_{y_o^i,o})=\partial_{y_o^i,p}$.

More precisely one has:
$$\tau_{\gamma}(\partial_{x_o^i,o})=\sum_j
\frac{\partial{x_r^j}}{\partial{x_o^i}}(p)\partial_{x_o^j,o},$$ et
$$\tau_{\gamma}(\partial_{y_o^i,o})=\sum_j
\frac{\partial{y_o^j}}{\partial{y_q^i}}(p)\partial_{y_o^j,o}.$$

Moreover the holonomy group is generated by the $\tau_\gamma$ of
this type.
\end{prop}
\begin{proof}
The first part of the proposition can be shown by evaluating
parallel transport of the given basis vectors along the given paths.

For the latter statement we will show show that one can approximate
any path $\gamma$ contained in a neighbourhood of $o$by a sequence
of paths $\gamma_n$ such that $\tau_{\gamma_n}$ is generated by the
parallel transports of curves of the type described in the
proposition (which we will call "rectangles" for simplicity) and
such that $\tau_\gamma$ is the limit of the $\tau_{\gamma_n}$.

Let $p_{1,p}$ be the local projection of $\mathcal M$ onto the leaf
${\mathcal F}_1(p)$ along the leaves of ${\mathcal F}_2$ and
$p_{2,p}$ be the local projection of $\mathcal M$ onto the leaf
${\mathcal F}_2(p)$ along the leaves of ${\mathcal F}_1$.

For a given curve $\mu:[t_0,t_1]\to M$ define the curve
$\mu_{(1)}:[t_0,t_1]\to M$ by $\mu_{(1)}=(p_{1,\mu(t_0)} \circ \mu)$
and the curve $\mu_{(2)}:[t_0,t_1]\to M$ by
$\mu_{(2)}=(p_{2,\mu(t_1)} \circ \mu)$. Finally let
$\mu':[t_0,t_1]\to M$ be the curve defined by $\mu'(\xi t_0 +
(1-\xi) \frac{t_0+t_1}{2})=\mu_{(1)}(\xi t_0 + (1-\xi) t_1)$ for $0
\le \xi \le 1$ and by $\mu'(\xi \frac{t_0+t_1}{2} +
(1-\xi)t_1)=\mu_{(2)}(\xi t_0 + (1-\xi) t_1)$ for $0\le \xi \le 1$.

For a given sequence of reals $0=t_0 < t_1 < t_2 < \ldots < t_n <
t_{n+1}=1$ one considers the curve $\nu:[0,1]\to M$ such that
$\nu\vrule_{[t_i,t_{i+1}]}=(\gamma\vrule_{[t_i,t_{i+1}]})'$ is
obtained by the preceding method.

One can show that the parallel transport $\tau_\nu$ is equal to a
product of parallel transports along "rectangles": In fact if
$\rho=\rho_1\rho_2$ is a closed curve then  $\rho'_1\rho'_2$ is a
"rectangle". One decomposes parallel transport $\tau_\nu$ then into
a product of parallel transports along the "rectangles": $R_1^{-1},
R_2, R_3^{-1}, R_4 \ldots, R_{2(n-1)}^{-1}, R_{2n}$ defined as
follows:

By setting $\omega_{2i}= \gamma\vrule_{[t_0,t_i]}$ and
$\omega_{2i+1}= \gamma\vrule_{[t_0,t_i]}
(\gamma\vrule_{[t_i,t_{i+1}]})_{(1)}$ for $i<n$, let $R_{r}$ be the
"rectangle" $\omega_{r}'(\omega_{r}^{-1})'$ for $r<2n$ and let
$R_{2n}$ be the "rectangle"
$(\gamma\vrule_{[t_0,t_{n}]})'(\gamma\vrule_{[t_n,t_{n+1}]})'$.

When $| t_{i+1}- t_i |$ tends towards $0$, $\nu$ tends to $\gamma$
in the sense of the compact-open topology which implies that
$\tau_\nu$ tends towards $\tau_\gamma$.

So the parallel transports along "rectangles" generate the whole
holonomy group.
\end{proof}

\subsection{Generic calculus}

\begin{prop}\label{cal_gen}
If the holonomy representation of the manifold $\mathcal M$ equipped
with the connection $\nabla$ is of type $T_o M=V_1 \oplus V_2$ with
$V_1$ and $V_2$ two faithful representations of the holonomy algebra
$\g g$, using the notations of pre preceding subsection and writing
$a_i^j$ and ${\tilde a}_l^m$ the coefficients defined by the
relations $\partial_{y_q^k,p}=\sum_{l}
a_{k}^{l}(p)\partial_{y_o^l,p},$ and $\partial_{y_o^l,p}=\sum_m
{\tilde a}_l^m(p)
\partial_{y_q^m,p}$
where $q$ is the intersection point of ${\mathcal F}_2(p)$ with
${\mathcal F}_1(o)$, defining $\Gamma$ by the relations
$$\nabla_{\partial_{x_o^i}}{\partial_{x_o^k}} = \sum_m
\Gamma_{i,k}^m {\partial_{x_o^m}},$$ one has:
\begin{eqnarray*}
R_p(\partial_{x_o^i},\partial_{y_o^j})\partial_{x_o^k} & = & -
\sum_m
\partial_{y_o^j} \left(\Gamma_{i,k}^m\right)\partial_{x_o^m},\\
R_p(\partial_{x_o^i},\partial_{y_o^j})\partial_{y_q^k} & = & -
\sum_{m}
\partial_{y_o^j}  \left( \sum_{l} \partial_{x_o^i} a_{k}^{l} {\tilde
a}_l^m \right)\partial_{y_q^m}.
\end{eqnarray*}
\end{prop}
\begin{proof} One can write (using in the point $p$ the equality:
$dy_o^i=\sum_j (\partial_{y_q^j}y_o^i)dy_q^j$):
$a_{k}^{l}(p)=\partial_{y_q^k,p}y_o^l$.

Clearly one has: $\sum_j a_k^j \tilde a_j^l=\delta_k^l$ and $\sum_j
\tilde a_k^j a_j^l=\delta_k^l$. By applying a derivation $\partial$
to these relations one obtains:

$\sum_j (\partial a_k^j) \tilde a_j^l= - \sum_j a_k^j (\partial
\tilde a_j^l)$ et $\sum_j (\partial \tilde a_k^j) a_j^l= - \sum_j
\tilde a_k^j (\partial a_j^l)$.

The equality $\nabla_{\partial_{y_o^j}}\partial_{y_q^k}=0$ (which is
verified because the leaf ${\mathcal F}_2(q)$ is flat), can be
written when defining $\Gamma'$ by the relations
$$\nabla_{\partial_{y_o^i}}{\partial_{y_o^k}} = \sum_m
{\Gamma'}_{i,k}^m {\partial_{y_o^m}}:$$

\begin{eqnarray*}
0 & = & \nabla_{\partial_{y_o^j}}\partial_{y_q^k} \\
 & = & \nabla_{\partial_{y_o^j}}\sum_{l} a_{k}^{l} \partial_{y_o^l}\\
 & = & \sum_{l}\left( (\partial_{y_o^j}a_{k}^{l})\partial_{y_o^l} +
a_{k}^{l}\nabla_{\partial_{y_o^j}}\partial_{y_o^l}\right)\\
& = & \sum_{l}\left( (\partial_{y_o^j}a_{k}^{l})\partial_{y_o^l} +
a_{k}^{l}\sum_p {\Gamma'}_{j,l}^p \partial_{y_o^p} \right)
\end{eqnarray*}

By separating the basis vectors one has for any $j,k,l$:

$$\partial_{y_o^j}a_{k}^{l} + \sum_p
a_{k}^{p}{\Gamma'}_{j,p}^l=0.\;\; (A)$$

One can isolate ${\Gamma'}_{j,p}^l$ by writing:
$${\Gamma'}_{j,p}^l=-\sum_k \tilde a_p^k \partial_{y_o^j}a_{k}^{l}. \;\; (B)$$

Deriving relation $(A)$ by $\partial_{x_o^i}$ one obtains:

$$\partial_{x_o^i}\partial_{y_o^j}a_{k}^{l} + \sum_p
(\partial_{x_o^i}a_{k}^{p}){\Gamma'}_{j,p}^l+\sum_p
a_{k}^{p}(\partial_{x_o^i}{\Gamma'}_{j,p}^l)=0.\;\; (*)$$

By choosing a path $\gamma_1$ from $o$ to $q$ contained in
${\mathcal F}_2(o)$ and a path $\gamma_2$ from $q$ to $p$ contained
in ${\mathcal F}_1(p)$ and noting $\gamma=\gamma_1 \gamma_2$ and
$\tau_\gamma$ the parallel transport along the path $\gamma$,
$\tau_\gamma(\partial_{x_o^k})=\partial_{x_o^k}$ and
$\tau_\gamma(\partial_{y_o^l})=\partial_{y_q^l}$.

At point $p$ one can calculate:
$$R_p(\partial_{x_o^i},\partial_{y_o^j})\partial_{x_o^k}=
\sum_m (-\partial_{y_o^j} \Gamma_{i,k}^m)\partial_{x_o^m} .$$

The more one has:
\begin{eqnarray*}
R_p(\partial_{x_o^i},\partial_{y_o^j})\partial_{y_q^k} & = &
\nabla_{\partial_{x_o^i}}
\underbrace{\nabla_{\partial_{y_o^j}}\partial_{y_q^k}}_{0} -
\nabla_{\partial_{y_o^j}}
\nabla_{\partial_{x_o^i}}\partial_{y_q^k}\\
 & = & - \nabla_{\partial_{y_o^j}}
\nabla_{\partial_{x_o^i}}\sum_l
a_{k}^{l}\partial_{y_o^l}\\
 & = & - \nabla_{\partial_{y_o^j}} \sum_l
(\partial_{x_o^i} a_{k}^{l})\partial_{y_o^l}-\underbrace{\sum_l
a_{k}^{l} \nabla_{\partial_{x_o^i}}\partial_{y_o^l}}_{0}\\
 & = & - \sum_l \left( (\partial_{y_o^j}\partial_{x_o^i}
a_{k}^{l})\partial_{y_o^l}+ \partial_{x_o^i}
a_{k}^{l}\nabla_{\partial_{y_o^j}}\partial_{y_o^l}\right) \\
 & = & - \sum_l \left( (\partial_{y_o^j}\partial_{x_o^i}
a_{k}^{l})\partial_{y_o^l}+ \partial_{x_o^i} a_{k}^{l}\sum_p
{\Gamma'}_{j,l}^p
\partial_{y_o^p}\right)
\end{eqnarray*}

By the relation $(*)$ and by the fact that $\partial_{x_o^i}$ and
$\partial_{y_o^j}$ commute one obtains:

\begin{eqnarray*}
R_p(\partial_{x_o^i},\partial_{y_o^j})\partial_{y_q^k} & = & \sum_l
\left(\sum_p a_{k}^{p}(\partial_{x_o^i}{\Gamma'}_{j,p}^l)\right)
\partial_{y_o^l}\\
 & = & \sum_{m,l,p} \left(
a_{k}^{p}(\partial_{x_o^i}{\Gamma'}_{j,p}^l) {\tilde a}_l^m \right)
\partial_{y_q^m},
\end{eqnarray*}

The coefficient of $\partial_{y_q^m}$ can be written:
\begin{eqnarray*}
\sum_{p,l} \left( a_{k}^{p}(\partial_{x_o^i}{\Gamma'}_{j,p}^l)
{\tilde a}_l^m \right) & = & \sum_{p,l} \left( a_{k}^{p}
\partial_{x_o^i}
(- \sum_n \tilde a_p^n \partial_{y_o^j}a_{n}^{l}) {\tilde a}_l^m
\right)\\
& = & - \sum_{p,n,l}\left( a_{k}^{p} \partial_{x_o^i}\tilde a_p^n
\partial_{y_o^j}a_{n}^{l}{\tilde a}_l^m + a_{k}^{p} \tilde a_p^n \partial_{x_o^i} \partial_{y_o^j}a_{n}^{l} {\tilde
a}_l^m \right)\\
& = & - \sum_{p,n,l}\left( \partial_{x_o^i} a_{k}^{p} \tilde a_p^n
a_{n}^{l}\partial_{y_o^j}{\tilde a}_l^m + a_{k}^{p} \tilde a_p^n
\partial_{x_o^i}
\partial_{y_o^j}a_{n}^{l} {\tilde a}_l^m \right)\\
& = & - \sum_{l}\left( \partial_{x_o^i} a_{k}^{l}
\partial_{y_o^j}{\tilde a}_l^m +
\partial_{x_o^i}
\partial_{y_o^j}a_{k}^{l} {\tilde a}_l^m \right)\\
& = & - \partial_{y_o^j}  \left( \sum_{l} \partial_{x_o^i} a_{k}^{l}
{\tilde a}_l^m \right)
\end{eqnarray*}
\end{proof}

\section{$V\oplus V^*$ holonomy}

In the para-Kähler case there is a bilinear symmetric non degenerate
form $b$ on $T_o M$ which is invariant by the action on the holonomy
algebra and $T_o M$ is equal to the direct sum of two totally
isotropic invariant subspaces $V_1$ and $V_2$. In particular the
signature of $b$ is neutral.

Remark that the mapping $\Psi: V_2 \to V_1^*$, $v\mapsto
b(v,\cdot)\vrule_{V_1}$ is injective, and surjective for dimension
reasons. In addition $\Psi$ is a morphism of representations. As a
consequence we have: $T_o M=V_1 \oplus V_1^*$.

Inversely if one has the equality of representations $E=V_1 \oplus
V_1^*$, one can associate to $E$ the symmetric bilinear non
degenerate form $b$ defined for $x,y\in V_1, u,v\in V_1^*$, by
$b(x+u,y+v)=u(y)+v(x)$. $V_1$ and $V_1^*$ are totally isotropic for
this form. As a conclusion we are in the para-Kähler case.

So it is equivalent to be in the para-Kähler case or to have a
holonomy of type $V\oplus V^*$.

Note that is this case the form $\omega(x+u,y+v)=u(y)-v(x)$ for
$x,y\in V_1, u,v\in V_1^*$ which is non degenerate bilinear
antisymmetric is also invariant.

\subsection{Berger criteria}
\begin{prop}\index{critère de Berger}
When $E=V\oplus V^*$, note $\tilde k(\g g) := (S^2(V^*)\otimes
S^2(V)) \cap (V^* \otimes \g g \otimes V)$, and $\underline{\g g} :=
<\{r(x,\cdot,\cdot,z') \; \vrule \; r\in \tilde k(\g g), x\in V,
z'\in V^* \}>$. If $\g g \subset \g g \g l(V)$ acting on $V \oplus
V^*$ is a holonomy representation of a torsion-free connection then
$\g g =\underline{\g g}$.
\end{prop}
\begin{proof} Recall that for $x,y \in V$, $z',t'\in V^*$ and $R$ a formal curvature tensor one has:
$R(x,y)=0$ and $R(z',t')=0$. $R$ is entirely known by the data of
$R(x,z')y$ for $x,y \in V$ and $z'\in V^*$. By the first Bianchi
identity necessarily $x,y \in V$ and $z'\in V^*$,
$R(x,z')y=R(y,z')x$. We can restate this by saying:
$R(\cdot,z')\cdot\vrule_{V\times V} \in \g g^{(1)}$.

Let us introduce the tensor $r$ defined for $x,y,z',t'\in V \oplus
V^*$ by $r(x,y,z',t'):=\langle R(x,z')y,t' \rangle$ which satisfies
the classical relations in pseudo-Riemannian geometry:
$r(x,y,z',t')=-r(z',y,x,t')=r(y,x,t',z')$. We can restrict to
$x,y\in V$ and $z',t'\in V^*$. As a consequence $r\in
S^2(V^*)\otimes S^2(V)$. Because$r(x,\cdot,\cdot,z')=R(x,z')\in \g
g$ it is clear that $r\in \tilde k(\g g)$.

A condition for $\g g$ being a holonomy representation is that $\g
g$ is generated by the $R(x,z')=r(x,\cdot,\cdot,z')$ for $r$ in
$\tilde k(\g g)$, $x$ in $V$ and $z'$ in $V^*$.
\end{proof}
\begin{prop}\index{critère de Berger}
When $E=V\oplus V^*$, note $$\tilde k^1(\g g) :=(S^3(V^*)\otimes
S^2(V)) \cap (S^2(V^*) \otimes \g g \otimes V) \oplus
(S^2(V^*)\otimes S^3(V)) \cap (V^* \otimes \g g \otimes S^2(V)).$$

If $\g g\subset \g g \g l(V)$ acting on $V \oplus V^*$ is a holonomy
representation of a torsion-free non locally symmetric connection
then $\tilde k^1(\g g)\neq 0$.
\end{prop}
\begin{proof}
This is again simply a specialization of the second Berger criterion
to this case.
\end{proof}

\subsection{Classification in dimension $2 \times 1$ and $2\times 2$}

In dimension $1+1$ the algebra $\g s \g o(1,1)$ corresponding to the
metric of matrix $\left(
\begin{smallmatrix}
0 & 1\\
1 & 0
\end{smallmatrix} \right)$ is precisely $\{ \left(
\begin{smallmatrix}
a & 0\\
0 & -a
\end{smallmatrix} \right) \; \vrule \; a\in \reels \}$.
It is the generic holonomy algebra in the pseudo-Riemannian case of
signature $(1,1)$.

The paper~\cite{LBBAI} gives a classification in signature $(2,2)$:

\begin{prop}
For $V$ of dimension $2$, $V\oplus V^*$ is an indecomposable
holonomy representation if and only if $V$ or $V^*$ is in the
following list:
\begin{enumerate}
\item $\left\{ \left(
\begin{smallmatrix}
0 & c\\
0 & 0
\end{smallmatrix} \right) \; \vrule \; c\in\reels\right\}$

\item For fixed $\lambda\in[-1,1]$, $\left\{ \left(
\begin{smallmatrix}
a & c\\
0 & \lambda a
\end{smallmatrix} \right) \; \vrule \; a,c\in\reels\right\}$

\item $\left\{ \left(
\begin{smallmatrix}
a & c\\
0 & b
\end{smallmatrix} \right) \; \vrule \; a,b,c\in\reels\right\}$

\item $\g s \g l(2,\reels)=\left\{ \left(
\begin{smallmatrix}
a & b\\
c & -a
\end{smallmatrix} \right) \; \vrule \; a,b,c\in\reels\right\}$

\item $\g g \g l(2,\reels)=\left\{ \left(
\begin{smallmatrix}
a & b\\
c & d
\end{smallmatrix} \right) \; \vrule \; a,b,c,d\in\reels\right\}$

\item $\g c \g o(2)=\left\{ \left(
\begin{smallmatrix}
a & -b\\
b & a
\end{smallmatrix} \right) \; \vrule \; a,b\in\reels\right\}$

\end{enumerate}

\end{prop}

\section{$V\oplus V$ holonomy}
A representation $V\oplus V$ can be written as well $V\otimes\reels
\oplus V\otimes\reels =V\otimes\reels^2$ where $\reels$ and
$\reels^2$ are the trivial representations.

\subsection{Berger criteria}
Let $V_1 \oplus V_2$ be a holonomy representation of $\g g \subset
\g g \g l(V_1)$ and let $\Phi:V_1 \to V_2, x \mapsto x'$ be an
isomorphism of representations. The curvature tensor $R$ verifies
then the following relations due to the Bianchi identity: $x,y,z\in
V_1:$
$$R(x,y)=0,$$
$$R(x',y')=0,$$
$$R(x,y')z=R(z,y')x,$$
$$R(x,y')z'=R(x,z')y'.$$
In addition because of the isomorphism $\Phi$ we have:
$$R(x,y')z'=(R(x,y')z)'.$$
One can deduce that the tensor $T$ defined for $x,y,z\in V_1$ by
$T(x,y,z):=R(x,y')z$ verifies:
$$T\in S^3(V^*_1)\otimes V_1.$$
In addition The Ambrose-Singer theorem gives us for $x,y\in V_1$:
$$T(x,y,\cdot)\in \g g.$$
Note $$\g g^{(2)}:= S^3(V^*_1)\otimes V_1 \cap (V^*_1)^{\otimes 2}
\otimes \g g$$ and $$\hat{\g g}:=\{T(x,y,\cdot) \; \vrule \; T\in \g
g^{(2)}, x, y \in V_1 \}.$$

The first Berger criterion can be formulated in this context simply
by
\begin{prop}
Let $\g g \subset \g g \g l(V_1)$. A necessary condition for $\g g
\subset \g g \g l(V_1)$ acting by $V_1 \otimes \reels^2$ is a
holonomy representation for a torsion-free connection is that $\g g
= \hat{\g g}$.
\end{prop}

Note that the second Bianchi identity can be formulated here by: for
$x,y,z\in V_1$,
$$(\nabla_x R)(y,z')=(\nabla_y R)(x,z'),$$ so that the tensor
$S$ defined by $S(x,y,z,t):=(\nabla_x R)(y,z')t$ lives in
$S^4(V^*_1)\otimes V_1$.

By noting $$\g g^{(3)}:= S^4(V^*_1)\otimes V_1 \cap (V^*_1)^{\otimes
3} \otimes \g g,$$ the second Berger criterion can be formulated:

\begin{prop}
If $\g g \subset \g g \g l(V_1)$ acting by $V_1 \otimes \reels^2$ is
a holonomy representation for a torsion-free non locally symmetric
connection then necessarily $\g g^{(3)}\neq \{0\}$.
\end{prop}

\subsection{Geometric structure}

By using the notations and calculus used in the proof of
proposition~\ref{cal_gen} we obtain:

If the holonomy is of type $V\otimes \reels^2$, at the origin $o$
$\partial_{x_o^k,o}$ and $\partial_{y_o^k,o}$ transform by the same
endomorphism. Their parallel transport along a path $\gamma$ from
$o$ to $p$ defined as the product $\gamma_1 \gamma_2$ of a path
$\gamma_1$ from $o$ to $q$ contained in ${\mathcal F}_1(o)$ and a
path $\gamma_2$ from $q$ to $p$ contained in ${\mathcal F}_2(p)$
gives: $\tau_\gamma(\partial_{x_o^k})=\partial_{x_o^k}$ et
$\tau_\gamma(\partial_{y_o^k})=\partial_{y_q^k}$.

$R_p(\partial_{x_o^i},\partial_{y_o^j})$ acts again "identically" on
these two vectors, which means that the relation between
$R_p(\partial_{x_o^i},\partial_{y_o^j})\partial_{x_o^k}$ and
$R_p(\partial_{x_o^i},\partial_{y_o^j})\partial_{y_q^k}$ obtained in
proposition~\ref{cal_gen} can be translated by:

$$\partial_{y_o^j} \Gamma_{i,k}^m = \partial_{y_o^j} (\sum_{l} \partial_{x_o^i} a_{k}^{l}
{\tilde a}_l^m),$$

or by: $$\Gamma_{i,k}^m = \sum_{l}
\partial_{x_o^i} a_{k}^{l} {\tilde a}_l^m + \mbox{const.}$$

By evaluating this equality in $o$ it follows
$$\Gamma_{i,k}^m = \sum_{l}
\partial_{x_o^i} a_{k}^{l} {\tilde a}_l^m.$$

So that we can write:
\begin{prop}
If the holonomy of the connection $\nabla$ is of type $V\otimes
\reels^2$, $a_i^j$ et ${\tilde a}_l^m$ being the coefficients
defined by $\partial_{y_q^k,p}=\sum_{l}
a_{k}^{l}(p)\partial_{y_o^l,p},$ and $\partial_{y_o^l,p}=\sum_m
{\tilde a}_l^m(p)
\partial_{y_q^m,p}$
(where $q$ is the intersection point of ${\mathcal F}_2(p)$ with
${\mathcal F}_1(o)$), defining $\Gamma$ by the relations
$$\nabla_{\partial_{x_o^i}}{\partial_{x_o^k}} = \sum_m
\Gamma_{i,k}^m {\partial_{x_o^m}},$$ we have:

$$\Gamma_{i,k}^m = \sum_{l}
\partial_{x_o^i} a_{k}^{l} {\tilde a}_l^m.$$
\end{prop}

The condition that $\nabla$ is torsion-free gives then the following
equalities:

\begin{prop}\label{p2}
The coordinates defined before verify the relations (for any $i,j,
k$):

$$\frac{\partial^2 y_o^k}{\partial x^i \partial y^j}=\frac{\partial^2 y_o^k}{\partial x^j \partial
y^i}.$$

\end{prop}

\subsection{Classification in dimension $2\times 2$}

Remark that the representations $V\otimes \reels^2$ which are
decomposable are all holonomy representations (for a torsion-free
connection) because all representations of dimension $2$ are.

\begin{prop}
For $V$ of dimension $2$, $V\otimes \reels^2$ is an indecomposable
holonomy representation if and only if $V$ is in the following list:
\begin{enumerate}
\item $\left\{ \left(
\begin{smallmatrix}
0 & c\\
0 & 0
\end{smallmatrix} \right) \; \vrule \; c\in\reels\right\}$

\item For fixed $\lambda\in \reels$, $\left\{ \left(
\begin{smallmatrix}
a & c\\
0 & \lambda a
\end{smallmatrix} \right) \; \vrule \; a,c\in\reels\right\}$

\item $\left\{ \left(
\begin{smallmatrix}
0 & c\\
0 & a
\end{smallmatrix} \right) \; \vrule \; a,c\in\reels\right\}$

\item $\left\{ \left(
\begin{smallmatrix}
a & c\\
0 & b
\end{smallmatrix} \right) \; \vrule \; a,b,c\in\reels\right\}$

\item $\g s \g l(2,\reels)=\left\{ \left(
\begin{smallmatrix}
a & b\\
c & -a
\end{smallmatrix} \right) \; \vrule \; a,b,c\in\reels\right\}$

\item $\g g \g l(2,\reels)=\left\{ \left(
\begin{smallmatrix}
a & b\\
c & d
\end{smallmatrix} \right) \; \vrule \; a,b,c,d\in\reels\right\}$

\item $\g c \g o(2)=\left\{ \left(
\begin{smallmatrix}
a & -b\\
b & a
\end{smallmatrix} \right) \; \vrule \; a,b\in\reels\right\}$

\end{enumerate}

\end{prop}

\begin{proof}
Recall the list of the representations of dimension $2$:

{\scriptsize
\begin{tabular}{|lll|ll|ll|}

\noalign{\hrule}
 & && Lie group && Lie algebra &\\
\noalign{\hrule}
  & $Gl^+(2,\reels)$\index{$Gl^+(2,\reels)$}
      &&
    $\left( \begin{smallmatrix}
    a & b \\
    c & d \end{smallmatrix}\right)$
    with $ad-bc>0$
      &&
    $\left( \begin{smallmatrix}
    a & b \\
    c & d \end{smallmatrix}\right)$
    with $a,b,c,d\in \reels$
      &\\
  \noalign{\hrule}

  & $Sl(2,\reels)$\index{$Sl(2,\reels)$}
      &&
    $\left( \begin{smallmatrix}
    a & b \\
    c & d \end{smallmatrix}\right)$
    with $ad-bc=1$
      &&
    $\left( \begin{smallmatrix}
    a & b \\
    c & -a \end{smallmatrix}\right)$
    with $a,b,c\in \reels$
      &\\
  \noalign{\hrule}

  & $SO(2)$\index{$SO(2)$}
      &&
    $\left( \begin{smallmatrix}
    \cos \theta & -\sin \theta \\
    \sin \theta & \cos \theta \end{smallmatrix}\right)$
    with $\theta \in \reels$
      &&
    $\left( \begin{smallmatrix}
    0 & a \\
    -a & 0 \end{smallmatrix}\right)$
    with $a\in \reels$
      &\\
  \noalign{\hrule}

  & $CO(2)$\index{$CO(2)$}
      &&
    $\left( \begin{smallmatrix}
    a & -b \\
    b & a \end{smallmatrix}\right)$
    with $a^2+b^2\neq 0$
      &&
    $\left( \begin{smallmatrix}
    a & -b \\
    b & a \end{smallmatrix}\right)$
    with $a,b\in \reels$
      &\\
  \noalign{\hrule}

  & $CO(2)_\lambda$\index{$CO(2)_\lambda$}
      &&
    $e^\theta\left( \begin{smallmatrix}
    \cos(\lambda \theta) & -\sin(\lambda \theta) \\
    \sin(\lambda \theta) & \cos(\lambda \theta) \end{smallmatrix}\right)$
    with $\lambda \neq 0$ fixed and $\theta\in\reels$
      &&
     $\left( \begin{smallmatrix}
    \theta & -\lambda \theta \\
    \lambda \theta & \theta \end{smallmatrix}\right)$
    avec $\theta\in \reels$
      &\\
  \noalign{\hrule}

  & $Id$\index{$Id$}
      &&
    $\left( \begin{smallmatrix}
    1 & 0 \\
    0 & 1 \end{smallmatrix}\right)$
      &&
    \quad $0$
      &\\
  \noalign{\hrule}

 & Homotheties\index{homothétie}
      &&
    $\left( \begin{smallmatrix}
    a & 0 \\
    0 & a \end{smallmatrix}\right)$
    with $a>0$
      &&
    $\left( \begin{smallmatrix}
    a & 0 \\
    0 & a \end{smallmatrix}\right)$
    with $a\in \reels$
      &\\
  \noalign{\hrule}

  & $SO(1,1)$\index{$SO(1,1)$}
      &&
    $\left( \begin{smallmatrix}
    a & 0 \\
    0 & a^{-1} \end{smallmatrix}\right)$
    with $a>0$
      &&
    $\left( \begin{smallmatrix}
    a & 0 \\
    0 & -a \end{smallmatrix}\right)$
    with $a\in \reels$
      &\\
  \noalign{\hrule}

  & $SO_\lambda$\index{$SO_\lambda$}
      &&
    $\left( \begin{smallmatrix}
    a & 0 \\
    0 & a^\lambda \end{smallmatrix}\right)$
    with $\lambda$ fixed $\neq -1,0,1$ and $a>0$
      &&
    $\left( \begin{smallmatrix}
    a & 0 \\
    0 & \lambda a \end{smallmatrix}\right)$
    with $a\in \reels$
      &\\
  \noalign{\hrule}

 & $SO-1_1$\index{$SO-1_1$}
      &&
    $\left( \begin{smallmatrix}
    a & 0 \\
    0 & 1 \end{smallmatrix}\right)$
    with $a>0$
      &&
    $\left( \begin{smallmatrix}
    a & 0 \\
    0 & 0 \end{smallmatrix}\right)$
    with $a\in \reels$
      &\\
  \noalign{\hrule}

  & $CO(1,1)$\index{$CO(1,1)$}
      &&
    $\left( \begin{smallmatrix}
    a & 0 \\
    0 & b \end{smallmatrix}\right)$
    with $a,b>0$
      &&
    $\left( \begin{smallmatrix}
    a & 0 \\
    0 & b \end{smallmatrix}\right)$
    with $a,b\in \reels$
      &\\
  \noalign{\hrule}

  & $He$\index{$He$}
      &&
    $\left( \begin{smallmatrix}
    1 & a \\
    0 & 1 \end{smallmatrix}\right)$
    with $a\in\reels$
      &&
    $\left( \begin{smallmatrix}
    0 & a \\
    0 & 0 \end{smallmatrix}\right)$
    with $a\in \reels$
      &\\
  \noalign{\hrule}

 & $Tr-X$\index{$Tr-X$}
      &&
    $\left( \begin{smallmatrix}
    e^t & t e^t \\
    0 & e^t \end{smallmatrix}\right)$
    with $t\in\reels$
      &&
    $\left( \begin{smallmatrix}
    a & a \\
    0 & a \end{smallmatrix}\right)$
    with $a\in \reels$
      &\\
  \noalign{\hrule}

  & $Tr-H$\index{$Tr-H$}
      &&
    $\left( \begin{smallmatrix}
    a & b \\
    0 & a\end{smallmatrix}\right)$
    with $a>0,b\in\reels$
      &&
    $\left( \begin{smallmatrix}
    a & b \\
    0 & a \end{smallmatrix}\right)$
    with $a,b\in \reels$
      &\\
  \noalign{\hrule}

  & $Tr-SO(1,1)$\index{$Tr-SO(1,1)$}
      &&
    $\left( \begin{smallmatrix}
    a & b \\
    0 & a^{-1}\end{smallmatrix}\right)$
    with $a>0,b\in\reels$
      &&
    $\left( \begin{smallmatrix}
    a & b \\
    0 & -a \end{smallmatrix}\right)$
    with $a,b\in \reels$
      &\\
  \noalign{\hrule}

  & $Tr-SO_\lambda$\index{$Tr-SO_\lambda$}
      &&
    $\left( \begin{smallmatrix}
    a & b \\
    0 & a^\lambda \end{smallmatrix}\right)$
    with $\lambda$ fixed $\neq -1,0,1$ and $a>0,b\in\reels$
      &&
    $\left( \begin{smallmatrix}
    a & b \\
    0 & \lambda a \end{smallmatrix}\right)$
    with $a,b\in \reels$
      &\\
  \noalign{\hrule}

  & $Tr-SO-1_1$\index{$Tr-SO-1_1$}
      &&
    $\left( \begin{smallmatrix}
    a & b \\
    0 & 1 \end{smallmatrix}\right)$
    with $a>0,b\in\reels$
      &&
    $\left( \begin{smallmatrix}
    a & b \\
    0 & 0 \end{smallmatrix}\right)$
    with $a,b\in \reels$
      &\\
  \noalign{\hrule}

  & $Tr-SO-1_2$\index{$Tr-SO-1_2$}
      &&
    $\left( \begin{smallmatrix}
    1 & b \\
    0 & a \end{smallmatrix}\right)$
    with $a>0,b\in\reels$
      &&
    $\left( \begin{smallmatrix}
    0 & b \\
    0 & a \end{smallmatrix}\right)$
    with $a,b\in \reels$
      &\\
  \noalign{\hrule}

  & $Tr$\index{$Tr$}
      &&
    $\left( \begin{smallmatrix}
    a & b \\
    0 & c \end{smallmatrix}\right)$
    with $a>0,c>0,b\in\reels$
      &&
    $\left( \begin{smallmatrix}
    a & b \\
    0 & c \end{smallmatrix}\right)$
    with $a,b,c\in \reels$
      &\\
  \noalign{\hrule}

\end{tabular}
}

By examining the different representations $V$ in dimension $2$, we
have:

$CO(1,1)$, $Id$ et $SO-1_1$ are decomposable.

$SO(2)$, $CO(2)_\lambda$ (for $\lambda\neq 0$), $SO(1,1)$, $Tr-X$,
the homotheties and $SO_\lambda$ (for $\lambda\neq -1,0,1)$) are
excluded because $V\otimes \reels^2$ is then a representation of a
Lie algebra of dimension $1$ whose image is generated by an
endomorphism of rank bigger or equal to $3$ which can not be a
holonomy representation because of the first Bianchi identity.

For the other cases we construct associated connections by
making explicit the coordinate relations by proposition~\ref{p1} and
verifying the relations stated in~\ref{p2}. Note that $y^i$ and
$y_o^i$ have to coincide if $\forall j, x^j=0.$

For $Gl^+(2,\reels)$ we take: $$y_o^1=x_o^1y^2+x_o^2y^1+y^1,$$
$$y_o^2=x_o^1y^2x_o^2y^1+ y^2.$$

$Sl(2,\reels)$ is an example known from the work of R. Bryant
(see~\cite{Br1} and \cite{Br2}): $\g s \g u(1,1)$ is a holonomy
representation.

For $CO(2)$ we take: $X=x_o^1+i x_o^2$, $Y=y^1+i y^2$ a holomorphic
function $F:U\to \complexes$ (with $U$ an open neighbourhood of $o$
in $\complexes^2$) such that $\frac{\partial^2 F}{\partial X
\partial Y}$ is not identically vanishing and such that
$F(0,Y)=Y$ et $F(X,0)=X$, for example $F(X,Y)=XY+X+Y$. We note
$F=F_R + i F_I$ with $F_R$ and $F_I$ real-valued functions. The
functions
$$y_o^1=F_R,$$
$$y_o^2=F_I$$
define then a connection of the wanted type.

$He$ is as well of type $V\oplus V^*$. One can take:
$$y_o^1=y^1+f(x^2,y^2),$$
$$y_o^2=y^2,$$
$f$ is a function vanishing for $x^2=0$.

For $Tr$ one can take

$$y_o^1=x^1y^1+y^1+y^2,$$
$$y_o^2=y^2(1+x^2),$$
as an example.

$Tr-H$, $Tr-SO(1,1)$, $Tr-SO_\lambda$ ($\lambda\neq -1,0,1$),
$Tr-SO-1_1$,which enter the scheme $\left(\begin{smallmatrix}a &
b\\0 & a^\lambda\end{smallmatrix}\right)$ with $\lambda$ one can
take:

$$y_o^1=x^1y^2+x^2y^1+y^1,$$
$$y_o^2=y^2(1+x^2)^\lambda.$$

For the group $Tr-SO-1_2$ one can take
$$y_o^1=y^1+y^2,$$
$$y_o^2=y^2(1+x^2).$$

\end{proof}

\section{Final remark}
We evaluated the Berger criteria for both types of representations
$V\oplus V^*$ and $V\otimes\reels^2$ for all representations $V$ up
to dimension $3$ using the software Maple. These results might
appear in future articles.



\begin{thebibliography}{CD9999}
\bibitem[AS]{AS}W. Ambrose, I. M. Singer, \textsl{A theorem on holonomy}, Trans. Amer. Math.
Soc. \textbf{79} (1953), pp. 428-443. \textsl{On the Holonomy of
Lorentzian Manifolds}, Proceeding of symposia in pure math., volume
54, 1993, pp. 27-40.
\bibitem[BBI]{LBBAI} L. Bérard Bergery, A. Ikemakhen,
\textsl{Sur l'holonomie des variétés pseudo-riemanniennes de
signature $(n,n)$}, Bulletin de la Société Mathématique de France,
125 no. 1 (1997), pp. 93-114.
\bibitem[BBK]{LBBTK} L. Bérard Bergery, T. Krantz, \textsl{Representations admitting two pairs of supplementary invariant
spaces}, arXiv:0704.2777/, 2007
\bibitem[Ber]{berger1} M. Berger, \textsl{Sur les groupes d'holonomie des variétés à connexion affine et des variétés riemanniennes}, Bull. Soc. Math. France \textbf{83} (1955), pp. 279-330.
\bibitem[BL]{bl} A. Borel, A. Lichnerowicz, \textsl{Groupes d'holonomie des variétés riemanniennes}, C.R. Acad. Sci. Paris \textbf{234} (1952), pp.1835-1837.
\bibitem[Boub]{Boubel} Ch. Boubel, \textsl{On the holonomy of Lorentzian metrics}, Prépublication UMPA no.
323 (2004), Écol. Norm. Sup. Lyon.
\bibitem[Br1]{Br1} R. L. Bryant, \textsl{Metrics with exceptional holonomy}, Ann. of Math. (2) \textbf{126} (1987), pp.
525-576, \textbf{MR} 89b:53084, \textbf{Zbl} 0637.53042
\bibitem[Br2]{Br2} R. L. Bryant, \textsl{Classical, exceptional, and exotic holonomies: a status report},
in Actes de la Table Ronde de Géométrie Différentielle (Luminy,
1992), Sémin. Congr., vol. 1 (1996), pp. 93–165, Soc. Math. France,
Paris \textbf{MR} 98c:53037. bibitem[Ca]{cartan} E. Cartan,
\textsl{Sur les variétés à connexion affine et la théorie de la
relativité générale I, II}, Ann. Sci. Écol. Norm. Sup. \textbf{40}
(1923), pp. 325-417, et \textbf{41} (1924), pp.1-25, Voir aussi
Œuvres complètes Tome III.
\bibitem[CMS]{cms}Q.-S. Chi, S. Merkulov, L. Schwachhöfer, \textsl{On the existence of infinite
series of exotic holonomies}, Inventiones Math. \textbf{126} (1996),
pp. 391–411.
\bibitem[DR]{DeRham} G. de Rham, \textsl{Sur la réductibilité d'un espace de Riemann}, Math Helv. \textbf{26} (1952), pp. 328-344.
\bibitem[Gal1]{Galaev1} A. Galaev, \textsl{Metrics that realize all
types of Lorentzian holonomy algebras}, arXiv:math.DG/0502575, 2005.
\bibitem[Gal2]{Galaev2} A. Galaev, \textsl{Holonomy groups and special geometric structures of pseudo-Kählerian manifolds
of index 2}, thèse de doctorat de la Humboldt-Universität, Berlin,
2006.
\bibitem[I2]{Ikemakhen2} A. Ikemakhen, \textsl{Examples of indecomposable non-irreducible Lorentzian
manifolds}, Ann. Sci. Math. Québec \textbf{20} (1996), no. 1, pp.
53-66.
\bibitem[I3]{Ikemakhen3}A. Ikemakhen, Sur l'holonomie des
variétés pseudo-riemanniennes de signature (2,2+n), 1998.
\bibitem[K]{krantz} T. Krantz, \textsl{Holonomie des connexions sans torsion}, PhD thesis of the Université Henri Poincaré Nancy I,
2007.
\bibitem[KN]{Kobayashi-Nomizu}
S.~Kobayashi, K.~Nomizu, \textsl{Foundations of differential
geometry, Vol I/II}, Interscience, 1963.
\bibitem[Le1]{leistner1} T. Leistner, \textsl{PhD Thesis, Humboldt University Berlin}, 2003.
\bibitem[Le2]{leistner2} T. Leistner, \textsl{Berger algebras, weak-Berger algebras and Lorentzian
holonomy}, sfb 288 preprint no. 567, 2002.
\bibitem[Le3]{leistner3} T. Leistner, \textsl{Towards a classification of Lorentzian holonomy groups}, arXiv:math.DG/0305139, 2003.
\bibitem[Le4]{leistner4} T. Leistner, \textsl{Towards a classification of Lorentzian holonomy groups. Part II:Semisimple, non-simple weak-Berger algebras}, arXiv:math.DG/0309274, 2003.
\bibitem[MS]{MS}S. Merkulov and L. Schwachhöfer, \textsl{Classification of irreducible
holonomies of torsion-free affine connections}, Ann. of Math. (2)
\textbf{150} (1999), no. 1, pp. 77-149.
\bibitem[Sch]{sch} L. Schwachhöfer, \textsl{On the classification of holonomy representations}, Habilitation Leipzig, 1999.
\bibitem[W]{Wu} H. Wu, \textsl{Holonomy groups of indefinite metrics},
Pacific Journal of Mathematics, \textbf{20} (1967), pp. 351-392.

\end{thebibliography}
\end{document}